\documentclass[10pt]{article} \usepackage{amsfonts} \usepackage{amssymb}\usepackage{amsmath}\usepackage{amsfonts}
\usepackage[russian]{babel}
\usepackage[cp1251]{inputenc}
{\Large

\makeatletter
\renewcommand{\@makefnmark}{}
\makeatother
\topmargin=-5mm
\oddsidemargin=-5mm
\evensidemargin=5mm
\textwidth=170mm
\textheight=220mm
\begin{document}
\baselineskip=20pt
\pagestyle{plain}

\makeatletter
\renewcommand{\@makefnmark}{}
\makeatother

\footnote{

\hspace{2mm}Keywords: Dirac operator, irregular and degenerate boundary conditions, completeness of root function systems}

\newcommand{\om}{O(\frac{1}{\mu})}
\newcommand{\lp}{L_2(0,\pi)}
\newcommand{\lpp}{L_2(\Omega)}

\newcommand{\elx}{e^{i\lambda x}}
\newcommand{\elxx}{e^{-i\lambda x}}

\newcommand{\elp}{e^{i\pi\lambda }}
\newcommand{\elpp}{e^{-i\pi\lambda }}
\newcommand{\No}{\textnumero}

\newcommand{\pp}{\prod_{n=-\infty\atop n\ne0}^\infty}
\newcommand{\sss}{\sum_{n=-\infty}^\infty}
\newcommand{\sN}{\sum_{n=-N}^N}
\newcommand{\sNN}{\sum_{|n|>N}}
\newcommand{\pps}{\sum_{n=-\infty\atop n\ne0}^\infty}
\newcommand{\sNNn}{\sum_{|n|\le N}}
\newcommand{\ppo}{\prod_{n=-\infty}^\infty}
\newcommand{\ppj}{\prod_{n=-\infty\atop n\ne0}^\infty\prod_{j=1}^2}
\newcommand{\ppjj}{\prod_{n=-\infty}^\infty\prod_{j=1}^2}

\newcommand{\ssj}{\sum_{n=-\infty\atop n\ne0}^\infty\sum_{j=1}^2}

\newcommand{\ssjj}{\sum_{n=-\infty}^\infty\sum_{j=1}^2}

\newcommand{\mun}{\mu_{n,j}}
\newcommand{\wnj}{W_{N,j}(\mu)}
\newcommand{\ntm}{|2n-\theta-\mu|}
\newcommand{\pwp}{PW_\pi^-}
\newcommand{\cpm}{\cos\pi\mu}
\newcommand{\spm}{\frac{\sin\pi\mu}{\mu}}
\newcommand{\tet}{(-1)^{\theta+1}}
\newcommand{\dm}{\Delta(\mu)}
\newcommand{\dl}{\Delta(\lambda)}

\newcommand{\dnm}{\Delta_N(\mu)}
\newcommand{\sni}{\sum_{n=N+1}^\infty}
\newcommand{\muk}{\sqrt{\mu^2+q_0}}
\newcommand{\agt}{\alpha, \gamma, \theta,}
\newcommand{\muq}{\sqrt{\mu^2-q_0}}
\newcommand{\smn}{\sum_{n=1}^\infty}
\newcommand{\lop}{{L_2(0,\pi)}}
\newcommand{\tdm}{\tilde\Delta_+(\mu_n)}
\newcommand{\dmn}{\Delta_+(\mu_n)}
\newcommand{\dd}{D_+(\mu_n)}
\newcommand{\ddd}{\tilde D_+(\mu_n)}
\newcommand{\sn}{\sum_{n=1}^\infty}
\newcommand{\pnn}{\prod_{n=1}^\infty}
\newcommand{\aln}{\alpha_n(\mu)}
\newcommand{\emp}{e^{\pi|Im\mu|}}
\newcommand{\ppp}{\prod_{p=p_0}^\infty}
\newcommand{\ppk}{\sum_{p=p_0}^{\infty}\sum_{k=1}^{[\ln p]}}
\newcommand{\mln}{m(\lambda_n)}
\newcommand{\mnk}{\mu_{n_k}}
\newcommand{\sxm}{s(x,\mu)}
\newcommand{\skm}{s(\xi,\mu)}
\newcommand{\cxm}{c(x,\mu)}
\newcommand{\ckm}{c(\xi,\mu)}

\newcommand{\lnk}{\lambda_{n_k}}
\medskip
\medskip
\medskip

\medskip
\medskip
\medskip
\medskip
\medskip

{\Large

\medskip
\medskip
\medskip
\medskip
\medskip

\centerline {\bf On the completeness of root function system of the Dirac operator }
\centerline {\bf with two-point boundary conditions }

\medskip
\medskip
\medskip
\medskip

\centerline { Alexander Makin}
\medskip
\medskip
\medskip

{\normalsize

\medskip
\medskip
\medskip
\medskip
\begin{quote}{{\bf Abstract.}
The paper is concerned with the completeness property of root functions of  the Dirac operator with summable complexvalued
 potential  and non-regular boundary conditions. We also obtain explicit form  for the fundamental solution system of the considered operator.}
  \end{quote}
\medskip
\medskip
\medskip
\medskip

\centerline {\bf 1. Introduction}
\medskip

In the  present paper, we study the Dirac system
\begin{equation}
B\mathbf{y}'+V\mathbf{y} =\lambda\mathbf{y},\label{1:2}
 \end{equation}
 where $\mathbf{y}={\rm col}(y_1(x),y_2(x))$,
 \[
 B=\begin{pmatrix}
 -i&0\\
 0&i
 \end{pmatrix},\quad V=\begin{pmatrix}
 0&P(x)\\
 Q(x)&0
 \end{pmatrix},
\]
the functions $P, Q\in L_1(0,\pi)$,   with two-point boundary conditions
$$
\begin{array}{c}
U_1(\mathbf{y})= a_{11}y_1(0)+a_{12}y_2(0)+a_{13}y_1(\pi)+ a_{14}y_2(\pi)=0,\\ U_2(\mathbf{y})= a_{21}y_1(0)+a_{22}y_2(0)+a_{23}y_1(\pi)+ a_{24}y_2(\pi)=0,
\end{array}\eqno(2)
$$
where
the coefficients $a_{jk}$ are arbitrary complex numbers,
and rows of  the matrix
\[
A=\begin{pmatrix}
 a_{11}&a_{12}&a_{13} &a_{14}\\
 a_{21}&a_{22}&a_{23} &a_{24}
\end{pmatrix}
\]
are linearly independent.

The operator $\mathbb{L}\mathbf{y}=B\mathbf{y}'+V\mathbf{y}$ is regarded as a linear operator in the space
$\mathbb{H}=L_2(0,\pi)\oplus L_2(0,\pi)$
with the domain $D(\mathbb{L})=\{\mathbf{y}\in W_1^1[0,\pi]\oplus  W_1^1[0,\pi]:\, \mathbb{L}\mathbf{y}\in \mathbb{H}$, $U_j(\mathbf{y})=0$ $(j=1,2)\}$.

Denote by $A_{jk}$ $(1\le j<k\le4)$ the determinant composed of the jth and kth columns of the matrix $A$.
It is known that boundary conditions (2) can be divided into three classes:

1)  regular conditions;

2)  irregular conditions;

3) degenerate conditions.

Boundary conditions (2) are called regular if
$$
A_{14}A_{23}\ne0.
$$
Boundary conditions (2) are called irregular if
$$
A_{14}A_{23}=0,\quad (A_{12}+A_{34})(|A_{23}|+|A_{14}|)\ne0
$$
and ones are called degenerate if
$$
A_{14}A_{23}=0,\quad (A_{12}+A_{34})(|A_{23}|+|A_{14}|)=0.
$$

It follows from [1] that
the root function system of problem (1), (2) with regular boundary  conditions is complete in $\mathbb{H}$.
If conditions (2) are not regular the completeness property essentially depends on the potential $V$, in particular, in this case
the root function system of nonperturbed operator

$$
B\mathbf{y}'=\lambda\mathbf{y},\quad U(\mathbf{y})=0\eqno(3)
$$
 is not complete in $\mathbb{H}$ [1].
In [2-5], A.V. Agibalova, A.A. Lunyov,  M.M. Malamud and L.L. Oridoroga  received a lot of results on the completeness (incompleteness) of root vectors for the Dirac operator (1), (2) and much more general first-order systems. In all mentioned papers the authors imposed
certain conditions on the smoothness of the coefficients of the considered operator. For example, for problem (1), (2) the functions $P, Q$ at least had to be continuous at the endpoints of the basic interval. Recently, in [6] A.P. Kosarev and A.A. Shkalikov established the completeness of the root function system of a problem of type (1), (2) if the coefficients are absolutely  continuous and satisfy some additional conditions at the endpoints of the basic interval.

 The aim of this paper is to examine the completeness property for the root function system of problem (1), (2) when the boundary  conditions are not regular and the potential $V\in L_1(0,\pi)$.

\medskip
\medskip

\centerline {\bf 2. Preliminaries}

Denote by
 \[
E(x,\lambda)=\begin{pmatrix}
e_{11}(x,\lambda)&e_{12}(x,\lambda)\\
e_{21}(x,\lambda)&e_{22}(x,\lambda)
\end{pmatrix}\eqno(4)
\]
the matrix of the fundamental solution system  to  equation
(1)
with boundary condition
$
E(0,\lambda)=I
$, where
$I$ is the unit matrix.
It is well known [7], [8] that
$$
\begin{array}{c}
e_{11}(x,\lambda)=e^{ix\lambda}(1+o(1))+e^{-i\pi\lambda}o(1),\quad e_{12}(x,\lambda)=e^{ix\lambda}o(1)+e^{-ix\lambda}o(1),\\

\quad e_{21}(x,\lambda)=e^{ix\lambda}o(1)+e^{-ix\lambda}o(1),\quad e_{22}(x,\lambda)=e^{ix\lambda}o(1)+e^{-ix\lambda}(1+o(1))
\end{array}\eqno(5)
$$
as $\lambda\to\infty$ uniformly in $x\in[0,\pi]$.

The eigenvalues of problem (1), (2) are the roots of the characteristic equation
$$
\Delta(\lambda)=0,
$$
where
$$
\Delta(\lambda)=
\left|\begin{array}{cccc}
U_1(E^{[1]}(\cdot,\lambda))&U_1(E^{[2]}(\cdot,\lambda))\\
U_2(E^{[1]}(\cdot,\lambda))&U_2(E^{[2]}(\cdot,\lambda))\\
\end{array}
\right|,
$$
$E^{[k]}(x,\lambda)$ is the $k$th column of matrix (4).

It was shown in [7] by the method of transformation operators that the characteristic determinant
$\Delta(\lambda)$ of problem (1), (2)
can be reduced to the form
$$
\begin{array}{c}
\Delta(\lambda)=A_{12}+A_{34}+A_{32}e_{11}(\pi,\lambda)+A_{14}e_{22}(\pi,\lambda)+A_{13}e_{12}(\pi,\lambda)+A_{42}e_{21}(\pi,\lambda)=\\
=\Delta_0(\lambda)+\int_0^\pi r_1(t)e^{-i\lambda t}dt+\int_0^\pi r_2(t)e^{i\lambda t}dt,
\end{array}\eqno(6)
$$
where the function
$$
\Delta_0(\lambda)=A_{12}+A_{34}-A_{23}e^{i\pi\lambda}+A_{14}e^{-i\pi\lambda}
$$
is the characteristic determinant of problem (3) and the functions $r_j\in L_1(0,\pi)$, $j=1,2$.

It is easy to see that in the case of degenerate boundary conditions the characteristic equation $\Delta_0(\lambda)=0$ either has no roots or
$\Delta_0(\lambda)\equiv0$.

 For convenience, we present several commonly used relations. Let $\lambda$ be a complex number, $Im\lambda\ne0$.
It follows from [9, p. 446] that
$$
\begin{array}{c}
\int_0^\infty x^\rho e^{2i\lambda x}dx=\\
\int_0^\infty x^\rho e^{-2Im\lambda x}(\cos(2Re\lambda x)+i\sin(2Re\lambda x))dx=\\=
\frac{\Gamma(\rho+1)}{2^{\rho+1}(Re^2\lambda+Im^2\lambda)^{(\rho+1)/2}}e^{i(\rho+1)\arctg(\frac{Re\lambda}{Im\lambda})}
\end{array}\eqno(7)
$$
if $Im\lambda>0$ and
$$
\begin{array}{c}
\int_0^\infty x^\rho e^{-2i\lambda x}dx=\\
\int_0^\infty x^\rho e^{2Im\lambda x}(\cos(2Re\lambda x)+i\sin(2Re\lambda x))dx=\\=
\frac{\Gamma(\rho+1)}{2^{\rho+1}(Re^2\lambda+Im^2\lambda)^{(\rho+1)/2}}e^{i(\rho+1)\arctg(\frac{Re\lambda}{-Im\lambda})}
\end{array}\eqno(8)
$$
if $Im\lambda<0$, where in both cases $\rho>-1$.

Let $\tau(x)$ be a continuous function on the segment $[0,\pi]$. It follows from (7), (8) that for any $b\in[0,\pi]$
$$
|\int_0^b x^\rho e^{-2|Im\lambda| x}\tau(x)dx|\le\frac{c}{|Im\lambda|^{\rho+1}},\eqno(9)
$$
where $c$  not depending on $b$.

It is easy to see that for any $b>0$
$$
\begin{array}{c}
|\int_b^\infty x^\rho e^{-2|Im\lambda| x}e^{2iRe\lambda x}dx|\le\int_b^\infty x^\rho e^{-2|Im\lambda| x}dx=\\
=\frac{1}{|Im\lambda|^{\rho+1}}\int_{b|Im\lambda|}^\infty t^\rho e^{-2t}dt=\frac{o(1)}{|Im\lambda|^{\rho+1}}
\end{array}\eqno(10)
$$
as $|Im\mu|\to\infty$.

Let $\tau(x)$ be a continuous function on the segment $[0,\pi]$, and $\tau(0)=0$.
Then, we have
$$
\begin{array}{c}
|\int_0^\pi x^\rho e^{-2|Im\lambda| x}\tau(x)dx|\le \int_0^\pi x^\rho e^{-2|Im\lambda| x}|\tau(x)|dx=\\
=\int_0^{\pi/\sqrt{|Im\lambda|}} x^\rho e^{-2|Im\lambda| x}|\tau(x)|dx
+\int_{\pi/\sqrt{|Im\lambda|}}^{\pi} x^\rho e^{-2|Im\lambda| x}|\tau(x)|dx\le\\
\le \max_{0\le x\le\pi/\sqrt{|Im\lambda|}}|\tau(x)|\int_0^\infty x^\rho e^{-2|Im\lambda| x}dx+\\
+e^{-\pi\sqrt{|Im\lambda|}}\int_0^\pi x^\rho|\tau(x)|dx=\frac{o(1)}{|Im\lambda|^{\rho+1}}
\end{array}\eqno(11)
$$
as $|Im\mu|\to\infty$.

Let a function $\tau(x)\in L_1(0,\pi)$. Then, we get

$$
\begin{array}{c}
|\int_0^\pi x^\rho e^{-2|Im\lambda| x}\tau(x)dx|\le \int_0^\pi x^\rho e^{-2|Im\lambda| x}|\tau(x)|dx=\\
=\int_0^{\pi/\sqrt{|Im\lambda|}} x^\rho e^{-2|Im\lambda| x}|\tau(x)|dx
+\int_{\pi/\sqrt{|Im\lambda|}}^{\pi} x^\rho e^{-2|Im\lambda| x}|\tau(x)|dx\le\\
\le \max_{0\le x\le\pi/\sqrt{|Im\lambda|}}( x^\rho e^{-2|Im\lambda| x})\int_0^{\pi/\sqrt{|Im\lambda|}}|\tau(x)|dx+e^{-\pi\sqrt{|Im\lambda|}}\int_0^\pi x^\rho|\tau(x)|dx\le\\
\le\frac{c_1}{|Im\lambda|^{\rho}}\int_0^{\pi/\sqrt{|Im\lambda|}}|\tau(x)|dx+e^{-\pi\sqrt{|Im\lambda|}}\int_0^\pi x^\rho|\tau(x)|dx=\frac{o(1)}{|Im\lambda|^{\rho}}
\end{array}\eqno(12)
$$
as $|Im\mu|\to\infty$.

Simple computations show that if $\rho>0, \lambda>0$,   $\rho\le \pi\lambda$, then
$$
\max_{0\le x\le \pi}x^\rho e^{-\lambda x}=\frac{\rho^\rho}{\lambda^\rho}e^{-\rho}. \eqno(13)
$$

\medskip
\medskip
\medskip

\newpage
\centerline {\bf 3. Main results.}
\medskip
\medskip
\medskip
First of all, we obtain explicit form  for the entries of matrix (4). To do this rewrite system (1) in scalar form
$$
\left\{
\begin{array}{rcl}
-iy_1'+P(x)y_2=\lambda y_1\\
iy_2'+Q(x)y_1=\lambda y_2.\\
\end{array}
\right.\eqno(14)
$$

 Multiplying  the first equation of (14)
by ${\elxx}$
and multiplying  the second equation of (14)
by ${\elx}$ and integrating the resulting equations over segment $[0,t]$, we obtain equivalent to (1)
the following system of integral equations

$$
\left\{
\begin{array}{rcl}
y_1(t)=e^{i\lambda t}y_1(0)-ie^{i\lambda t}\int_0^t{\elxx}P(x)y_2(x)dx\\
y_2(t)=e^{-i\lambda t}y_2(0)+ie^{-i\lambda t}\int_0^t{\elx}Q(x)y_1(x)dx.
\end{array}
\right.\eqno(15)
$$
It follows from (15) that

$$
\left\{
\begin{array}{rcl}
e_{11}(t,\lambda)=e^{i\lambda t}-ie^{i\lambda t}\int_0^te^{-i\lambda t_1}P(t_1)e_{21}(t_1,\lambda)dt_1\\
e_{21}(t,\lambda)=ie^{-i\lambda t}\int_0^te^{i\lambda t_1}Q(t_1)e_{11}(t_1,\lambda)dt_1,
\end{array}
\right.\eqno(16)
$$
hence,
$$
\begin{array}{c}
e_{11}(t,\lambda)=e^{i\lambda t}+e^{i\lambda t}\int_0^te^{-2i\lambda t_1}P(t_1)dt_1\int_0^{t_1}e^{i\lambda t_2}Q(t_2)e_{11}(t_2,\lambda)dt_2
\end{array}\eqno(17)
$$
and
$$
\begin{array}{c}
e_{21}(t,\lambda)=ie^{-i\lambda t}\int_0^te^{i\lambda t_1}Q(t_1)dt_1(e^{i\lambda t_1}-ie^{i\lambda t_1}\int_0^{t_1}e^{-i\lambda t_2}P(t_2)e_{21}(t_2,\lambda)dt_2)=\\
=ie^{-i\lambda t}\int_0^te^{2i\lambda t_1}Q(t_1)dt_1+e^{-i\lambda t}\int_0^te^{2i\lambda t_1}Q(t_1)dt_1\int_0^{t_1}e^{-i\lambda t_2}P(t_2)e_{21}(t_2,\lambda)dt_2.
\end{array}\eqno(18)
$$

In the same way, we get

$$
\left\{
\begin{array}{rcl}
e_{12}(t,\lambda)=-ie^{i\lambda t}\int_0^te^{-i\lambda t_1}P(t_1)e_{22}(t_1,\lambda)dt_1\\
e_{22}(t,\lambda)=e^{-i\lambda t}+ie^{-i\lambda t}\int_0^te^{i\lambda t_1}Q(t_1)e_{12}(t_1,\lambda)dt_1,
\end{array}
\right.
\eqno(19)
$$
hence,
$$
\begin{array}{c}
e_{22}(t,\lambda)=e^{-i\lambda t}+e^{-i\lambda t}\int_0^te^{2i\lambda t_1}Q(t_1)dt_1\int_0^{t_1}e^{-i\lambda t_2}P(t_2)e_{22}(t_2,\lambda)dt_2
\end{array}\eqno(20)
$$
and
$$
\begin{array}{c}
e_{12}(t,\lambda)=-ie^{i\lambda t}\int_0^te^{-i\lambda t_1}P(t_1)dt_1(e^{-i\lambda t_1}+ie^{-i\lambda t_1}\int_0^{t_1}e^{i\lambda t_2}Q(t_2)e_{12}(t_2,\lambda)dt_2)=\\
=-ie^{i\lambda t}\int_0^te^{-2i\lambda t_1}P(t_1)dt_1+e^{i\lambda t}\int_0^te^{-2i\lambda t_1}P(t_1)dt_1\int_0^{t_1}e^{i\lambda t_2}Q(t_2)e_{12}(t_2,\lambda)dt_2.
\end{array}
\eqno(21)
$$

Denote
$$
p_0(t,\lambda)=\int_0^t e^{-2i\lambda t_1}P(t_1)dt_1,\eqno(22)
$$

$$
p_1(t,\lambda)=\int_0^t e^{-2i\lambda t_1}P(t_1)dt_1\int_0^{t_1}e^{2i\lambda t_2}Q(t_2)dt_2\int_0^{t_2} e^{-2i\lambda t_3}P(t_3)dt_3,\eqno(23)
$$
$$
p_n(t,\lambda)=\int_0^t e^{-2i\lambda t_1}P(t_1)dt_1\int_0^{t_1}e^{2i\lambda t_2}Q(t_2)dt_2\ldots\int_0^{t_{2n}} e^{-2i\lambda t_{2n+1}}P(t_{2n+1})dt_{2n+1}\eqno(24)
$$
and, analogously,
denote
$$
q_0(t,\lambda)=\int_0^t e^{2i\lambda t_1}Q(t_1)dt_1,\eqno(25)
$$

$$
q_1(t,\lambda)=\int_0^t e^{2i\lambda t_1}Q(t_1)dt_1\int_0^{t_1}e^{-2i\lambda t_2}P(t_2)dt_2\int_0^{t_2} e^{2i\lambda t_3}Q(t_3)dt_3,\eqno(26)
$$
$$
q_n(t,\lambda)=\int_0^t e^{2i\lambda t_1}Q(t_1)dt_1\int_0^{t_1}e^{-2i\lambda t_2}P(t_2)dt_2\ldots\int_0^{t_{2n}} e^{2i\lambda t_{2n+1}}Q(t_{2n+1})dt_{2n+1}.\eqno(27)
$$

Denote also $\tilde p_n(t,\lambda)=-ie^{i\lambda t}p_n(t,\lambda)$, $\tilde q_n(t,\lambda)=ie^{-i\lambda t}q_n(t,\lambda)$.

Denote
$$
g_0(t,\lambda)=1,\eqno(28)
$$

$$
g_1(t,\lambda)=\int_0^t e^{-2i\lambda t_1}P(t_1)dt_1\int_0^{t_1}e^{2i\lambda t_2}Q(t_2)dt_2,\eqno(29)
$$
$$
\begin{array}{c}
g_n(t,\lambda)=\int_0^t e^{-2i\lambda t_1}P(t_1)dt_1\int_0^{t_1}e^{2i\lambda t_2}Q(t_2)dt_2\ldots\\\ldots
\int_0^{t_{2n-2}}e^{-2i\lambda t_{2n-1}}P(t_{2n-1})dt_{2n-1}\int_0^{t_{2n-1}}e^{2i\lambda t_{2n}}Q(t_{2n})dt_{2n}
\end{array}\eqno(30)
$$
and, analogously,
denote
$$
h_0(t,\lambda)=1,\eqno(31)
$$

$$
h_1(t,\lambda)=\int_0^t e^{2i\lambda t_1}Q(t_1)dt_1\int_0^{t_1}e^{-2i\lambda t_2}P(t_2)dt_2,\eqno(32)
$$
$$
\begin{array}{c}
h_n(t,\lambda)=\int_0^t e^{2i\lambda t_1}Q(t_1)dt_1\int_0^{t_1}e^{-2i\lambda t_2}P(t_2)dt_2\ldots\\\ldots
\int_0^{t_{2n-2}} e^{2i\lambda t_{2n-1}}Q(t_{2n-1})dt_{2n-1}\int_0^{t_{2n-1}}e^{-2i\lambda t_{2n}}P(t_2)dt_{2n}.
\end{array}\eqno(33)
$$

Denote also $\tilde g_n(t,\lambda)=e^{i\lambda t}g_n(t,\lambda)$, $\tilde h_n(t,\lambda)=e^{-i\lambda t}h_n(t,\lambda)$.

{\bf Lemma 1.} {\it Let a function $f\in L_1(a,x_0)$. Then, for any $n=1,2,\ldots$
$$
I_n=\int_a^{x_0}f(x_1)dx_1\int_a^{x_1}f(x_2)dx_2\ldots\int_a^{x_{n-1}}f(x_n)dx_n=\frac{(\int_a^{x_0}f(x)dx)^n}{n!}.\eqno(34)
$$}
The proof is by induction on $n$. For $n=1$, there is nothing to prove. Suppose equality (34) is valid for a number $n$. Then
$$
I_{n+1}=\int_a^{x_0}f(x_1)dx_1\int_a^{x_1}f(x_2)dx_2\ldots\int_a^{x_{n}}f(x_{n+1})dx_n=
\int_a^{x_0}f(x_1)dx_1\left(\frac{(\int_a^{x_1}f(x)dx)^n}{n!}\right).
$$
Integrating by parts, we obtain
$$
\begin{array}{c}
I_{n+1}=\frac{(\int_a^{x_0}f(x)dx)^{n+1}}{n!}-\frac{1}{(n-1)!}\int_a^{x_0}f(x_1)dx_1\left(\int_a^{x_1}f(x)dx)\right)^n=\\
=\frac{(\int_a^{x_0}f(x)dx)^{n+1}}{n!}-nI_{n+1},
\end{array}
$$
hence,
$$
(n+1)I_{n+1}=\frac{(\int_a^{x_0}f(x)dx)^{n+1}}{n!},
$$
therefore,
$$
I_{n+1}=\frac{(\int_a^{x_0}f(x)dx)^{n+1}}{(n+1)!}.
$$

{\bf Remark 1.} The similar argument one can meet in a lot of papers, for example,  in [10, 11].

{\bf Lemma 2.} {\it The following representations are valid
$$
e_{12}(t,\lambda)=-ie^{i\lambda t}\sum_{n=0}^\infty p_n(t,\lambda),\eqno(35)
$$
$$
e_{21}(t,\lambda)=ie^{-i\lambda t}\sum_{n=0}^\infty q_n(t,\lambda),\eqno(36)
$$

$$
e_{11}(t,\lambda)=e^{i\lambda t}\sum_{n=0}^\infty g_n(t,\lambda),\eqno(37)
$$
$$
e_{22}(t,\lambda)=e^{-i\lambda t}\sum_{n=0}^\infty h_n(t,\lambda),\eqno(38)
$$
where the series in right-hand sides of (35-38) for any $\lambda$ converge uniformly and absolutely on the segment $[0,\pi]$. }

Proof. Fix an arbitrary complex $\lambda$. Denote $\Phi(t)=|e^{-2i\lambda t}P(t)|+ |e^{2i\lambda t}Q(t)|$, $M=\int_0^\pi \Phi(t)dt$. Then, by virtue of Lemma 1
$$
\begin{array}{c}
|p_n(t,\lambda)|\le\int_0^t |e^{-2i\lambda t_1}P(t_1)|dt_1\int_0^{t_1}|e^{2i\lambda t_2}Q(t_2)dt_2|\ldots\int_0^{t_{2n}} |e^{-2i\lambda t_{2n+1}}P(t_{2n+1})|dt_{2n+1}\le\\ \le\int_0^t \Phi(t_1)dt_1\int_0^{t_1}\Phi(t_2)dt_2\ldots\int_0^{t_{2n}}\Phi(t_{2n+1})|dt_{2n+1}=
\frac{(\int_0^t \Phi(y)dy)^{2n+1}}{(2n+1)!}\le\frac{M^{2n+1}}{(2n+1)!}.
\end{array}
$$
 It is easy to see that for $n=0,1,\ldots$
$$
e^{i\lambda t}\int_0^te^{-2i\lambda t_1}P(t_1)dt_1\int_0^{t_1}e^{i\lambda t_2}Q(t_2)\tilde p_n(t_2,\lambda)dt_2=\tilde p_{n+1}(t,\lambda),\eqno(39)
$$
hence, substituting the series in  right-hand side of (35) in (21), we obtain the valid identity.

 In the same way, one can prove equalities (36-38). This completes the proof of  Lemma 2.

It follows from (39) that
$$
\int_0^te^{-2i\lambda t_1}P(t_1)dt_1\int_0^{t_1}e^{2i\lambda t_2}Q(t_2) p_n(t_2,\lambda)dt_2= p_{n+1}(t,\lambda).\eqno(40)
$$

Changing the order of integration in (40), we obtain
$$
p_{n+1}(t,\lambda)=\int_0^te^{i\lambda t_2}Q(t_2) p_n(t_2,\lambda)dt_2
\int_{t_2}^{t}e^{-2i\lambda t_1}P(t_1)dt_1.\eqno(41)
$$

To estimate the functions $e_{jk}(\pi,\lambda)$ $(1\le j,k\le2)$ we need the following

{\bf Lemma 3.} {\it Let a function $f\in L_1(0,\pi)$, $0\le a<b\le\pi$. Then, uniformly in $a,b$
$$
\lim_{|\lambda|\to\infty}e^{-|Im\lambda| b}\int_a^b e^{i\lambda x}f(x)dx=0,
\quad \lim_{|\lambda|\to\infty}e^{|Im\lambda| a}\int_a^b e^{-i\lambda x}f(x)dx=0
$$
if $Im\lambda<0$; and
$$
\lim_{|\lambda|\to\infty}e^{Im\lambda a}\int_a^b e^{i\lambda x}f(x)dx=0,\quad
\lim_{|\lambda|\to\infty}e^{-Im\lambda b}\int_a^b e^{-i\lambda x}f(x)dx=0
$$
if $Im\lambda\ge0$.}

Lemma 3 is an insignificant generalization of the well known Riemann lemma [1, Lemma 1.3.1].

Let $0<\varepsilon<\pi/10$.
Denote by $\Omega^+_\varepsilon$ the domain  $\varepsilon\le \arg\lambda\le\pi-\varepsilon$, and by $\Omega^-_\varepsilon$ the domain
 $-\pi+\varepsilon\le \arg\lambda\le-\varepsilon$.

Denote also $\tilde R_{12}(t,\lambda)=\sum_{n=1}^\infty \tilde p_n(t,\lambda)$, $\tilde R_{21}(t,\lambda)=\sum_{n=1}^\infty \tilde q_n(t,\lambda)$.

{\bf Lemma 4.} {\it Suppose
$$
\lim_{h\to0}\frac{\int_{\pi-h}^{\pi}P(x)dx}{h^\rho}=\nu\ne0,\eqno(42)
$$
where $\rho>0$. Then, uniformly in the domain $\Omega_\varepsilon^+$
$$
\begin{array}{c}
\tilde p_{0}(\pi,\lambda)=\frac{-\nu\lambda e^{-i\pi\lambda}\Gamma(\rho+1)}{2^{\rho}(Re^2\lambda+Im^2\lambda)^{(\rho+1)/2}}e^{i(\rho+1)\arctg(\frac{Re\lambda}{Im\lambda})}
+\frac{e^{\pi|Im\lambda|}o(1)}{|Im\lambda|^\rho},
\end{array}\eqno(43)
$$

$$
\tilde R_{12}(\pi,\lambda)=\frac{e^{\pi|Im\lambda|}o(1)}{|Im\lambda|^\rho},\eqno(44)
$$
and

$$
e_{11}(\pi,\lambda)=\frac{e^{\pi|Im\lambda|}o(1)}{|Im\lambda|^\rho}\eqno(45)
$$

as $|Im\mu|\to\infty$.

}
Proof. Denote    $\hat P(x)=P(\pi-x)$. It follows from (42) that
$$
\lim_{h\to0}\frac{\int_0^h\hat P(x)dx}{h^\rho}=\nu,
$$
hence,
$$
\int_{0}^{h}\hat P(x)dx=\nu h^\rho+h^\rho\tau(h),
$$
 where $\tau(h)\to0$ as $h\to0$.
Integrating by parts, we have
$$
\begin{array}{c}
 \int_0^\pi e^{-2i\lambda t}P(t)dt=e^{-2i\pi\lambda}\int_0^\pi e^{2i\lambda s}\hat P(s)ds=\\=

 \int_0^\pi \hat P(y)dy-2i\lambda e^{-2i\pi\lambda}\int_0^\pi e^{2i\lambda s}(\nu s^\rho+s^\rho\tau(s))ds=\\=

 \int_0^\pi \hat P(y)dy-2i\nu\lambda e^{-2i\lambda\pi}\int_0^\infty x^\rho e^{2i\lambda x}dx+2i\nu\lambda e^{-2i\lambda\pi}\int_\pi^\infty x^\rho e^{2i\lambda x}dx-\\
 -2i\lambda e^{-2i\lambda\pi}\int_0^\infty x^\rho e^{2i\lambda x}\tau(x)dx.
\end{array}
$$
This together with (10), (11) implies

$$
\int_0^\pi e^{-2i\lambda t}P(t)dt=-2i\nu\lambda e^{-2i\lambda\pi}\int_0^\infty t^\rho e^{-2i\lambda  t}dt+\frac{e^{2\pi|Im\lambda|}o(1)}{|Im\lambda|^{\rho}}.
$$
Multiplying the last relation by $-ie^{i\lambda \pi}$ and invoking (7), we obtain (43).

Let us estimate the function $p_n(\pi,\lambda)$.
Denote

 $$
 \begin{array}{c}
 F_n(t_2,\lambda)=e^{2i\lambda t_2}Q(t_2)\int_0^{t_{2}} e^{-2i\lambda t_{3}}P(t_{3})dt_{3}\ldots\\\ldots\int_0^{t_{2n-1}}e^{2i\lambda t_{2n}}Q(t_{2n})dt_{2n}\int_0^{t_{2n}} e^{-2i\lambda t_{2n+1}}P(t_{2n+1})dt_{2n+1}.
 \end{array}\eqno(46)
 $$

Replacing  $t_1=\pi-x$ in (24) and integrating by parts, we get

$$
\begin{array}{c}
p_n(\pi,\lambda)=e^{-2i\lambda \pi}\int_0^\pi e^{2i\lambda x}\hat P(x)dx\int_0^{\pi-x}F_n(t_2,\lambda)dt_2=\\=
e^{-2i\lambda \pi}\int_0^\pi (e^{2i\lambda x}\int_0^{\pi-x}F_n(t_2)dt_2)
d\int_0^{x}\hat P(y)dy=\\=
e^{-2i\lambda \pi}[\left.(e^{2i\lambda x}\int_0^{\pi-x}F_n(t_2)dt_2)\int_0^{x}\hat P(y)dy\right|_0^\pi-
\int_0^\pi (\int_0^{x}\hat P(y)dy)(e^{2i\lambda x}\int_0^{\pi-x}F_n(t_2,\lambda)dt_2)'dx]=\\=
e^{-2i\lambda \pi}[-2i\lambda\int_0^\pi(\nu x^\rho+x^\rho\tau(x))e^{2i\lambda x}(\int_0^{\pi-x}F_n(t_2,\lambda)dt_2)dx+\\+
\int_0^\pi(\nu x^\rho+x^\rho\tau(x))e^{2i\lambda x}F_n(\pi-x,\lambda)dx].
\end{array}\eqno(47)
$$
Let us estimate the first addend in square brackets in right-hand side of (47).
It follows from Lemma 3 that
$$
\begin{array}{c}
|\int_0^{t_{2n-1}}e^{2i\lambda t_{2n}}Q(t_{2n})dt_{2n}\int_0^{t_{2n}} e^{-2i\lambda t_{2n+1}}P(t_{2n+1})dt_{2n+1}|\le\\\le
\int_0^{t_{2n-1}}|Q(t_{2n}|)dt_{2n}
\max_{0\le t_{2n}\le t_{2n-1}}|e^{2i\lambda t_{2n}}\int_0^{t_{2n}} e^{-2i\lambda t_{2n+1}}P(t_{2n+1})dt_{2n+1}|\le\\\le
\epsilon(\lambda)\int_0^{t_{2n-1}}|Q(t_{2n}|)dt_{2n},
\end{array}\eqno(48)
$$
where $\epsilon(\lambda)\to0$ as $\lambda\to\infty$.
Since  $t_j\ge t_{j+1}$ for any $j$, it follows from (46), (48), and Lemma 1 that
$$
\begin{array}{c}
|\int_0^{\pi-x}F_n(t_2,\lambda)dt_2|\le \epsilon(\lambda)\int_0^{\pi-x}|Q(t_2)|\int_0^{t_{2}}| P(t_{3})|dt_{3}\ldots \\\ldots|e^{2i\lambda(t_2-t_3+\ldots +t_{2n-2}-t_{2n-1})}|\int_0^{t_{2n-1}}|Q(t_{2n})|dt_{2n} \le
\frac{M_1^{2n-1}}{(2n-1)!}\epsilon(\lambda).
\end{array}\eqno(49)
$$

The last inequality together with (9) implies
$$
\begin{array}{c}
|2i\lambda\int_0^\pi(\nu x^\rho+x^\rho\tau(x))e^{2i\lambda x}(\int_0^{\pi-x}F_n(t_2)dt_2)dx|\le
\frac{c_1M_1^{2n-1}\epsilon(\lambda)}{(2n-1)!|Im\lambda|^{\rho}}.
\end{array}\eqno(50)
$$
Let us estimate the second addend in square brackets in right-hand side of (47).
Clearly,

$$
\begin{array}{c}
\int_0^\pi(\nu x^\rho+x^\rho\tau(x))e^{2i\lambda x}F_n(\pi-x)dx=
\int_0^\pi(\nu x^\rho+x^\rho\tau(x))e^{2i\lambda x}Q(\pi-x)\phi_n(x,\lambda)dx,
\end{array}
$$
where
$$
\begin{array}{c}
\phi_n(x,\lambda)=e^{2i\lambda(\pi-x)}\int_0^{\pi-x} e^{-2i\lambda t_{3}}P(t_{3})dt_{3}\ldots\\\ldots\int_0^{t_{2n-1}}e^{2i\lambda t_{2n}}Q(t_{2n})dt_{2n}\int_0^{t_{2n}} e^{-2i\lambda t_{2n+1}}P(t_{2n+1})dt_{2n+1}.
\end{array}
$$
It follows from Lemma 1 that
$$
\begin{array}{c}
|\phi_n(x,\lambda)|\le \int_0^{\pi-x} |P(t_{3})|dt_{3}\ldots\\\ldots\int_0^{t_{2n-1}}|Q(t_{2n})|dt_{2n}\int_0^{t_{2n}} |e^{2i(\pi-x-t_3+\ldots+t_{2n}-t_{2n+1})\lambda }||P(t_{2n+1})|dt_{2n+1}\le\frac{M^{2n-1}}{(2n-1)!}.
\end{array}
$$
This and (12) imply
$$
\int_0^\pi(\nu x^\rho+x^\rho\tau(x))e^{2i\lambda x}F_n(\pi-x)dx=\frac{M^{2n-1}}{(2n-1)!}\frac{o(1)}{|Im\lambda|^{\rho}}.\eqno(51)
$$
It follows from (47), (50), (51)  that
$$
\sum_{n=1}^\infty|p_{n}(\pi,\lambda)|=\frac{e^{2\pi|Im\lambda|}o(1)}{|Im\lambda|^{\rho}},\eqno(52)
$$
therefore, inequality (44) holds.

Let us prove relation (45).
Denote
$$
G_1(t_2,\lambda)=e^{2i\lambda t_2}Q(t_2),
$$

 $$
 G_n(t_2,\lambda)=e^{2i\lambda t_2}Q(t_2)\int_0^{t_{2}} e^{-2i\lambda t_{3}}P(t_{3})dt_{3}\ldots\\\ldots
\int_0^{t_{2n-2}}e^{-2i\lambda t_{2n-1}}P(t_{2n-1})dt_{2n-1}\int_0^{t_{2n-1}}e^{2i\lambda t_{2n}}Q(t_{2n})dt_{2n}
 $$ $(n>1)$.
 Reasoning as above, we obtain

$$
\begin{array}{c}
g_n(\pi,\lambda)=\int_0^\pi e^{-2i\lambda t_1}P(t_1)dt_1\int_0^{t_1}G_n(t_2)dt_2=\\=
e^{-2i\lambda \pi}[-2i\lambda\int_0^\pi(\nu x^\rho+x^\rho\tau(x))e^{2i\lambda x}(\int_0^{\pi-x}G_n(t_2,\lambda)dt_2)dx+\\+
\int_0^\pi(\nu x^\rho+x^\rho\tau(x))e^{2i\lambda x}G_n(\pi-x,\lambda)dx].
\end{array}\eqno(53)
$$
Let us estimate the first addend in square brackets in right-hand side of (53).
It follows from Lemma 3 that
$$
\begin{array}{c}
|\int_0^{t_{2n-1}}e^{2i\lambda t_{2n}}Q(t_{2n})dt_{2n}|=\epsilon(\lambda),
\end{array}\eqno(54)
$$
where $\epsilon(\lambda)\to0$ as $\lambda\to\infty$.
Since  $t_j\ge t_{j+1}$ for any $j$, it follows from (54) and Lemma 1 that
$$
\begin{array}{c}
|\int_0^{\pi-x}G_n(t_2,\lambda)dt_2|\le \epsilon(\lambda)\int_0^{\pi-x}|Q(t_2)|\int_0^{t_{2}}| P(t_{3})|dt_{3}\ldots \\\ldots\int_0^{t_{2n-2}}|e^{2i\lambda(t_2-t_3+\ldots +t_{2n-2}-t_{2n-1})}||P(t_{2n-1})|dt_{2n-1} \le
\frac{M_2^{2n-1}}{(2n-1)!}\epsilon(\lambda).
\end{array}\eqno(55)
$$
Let us estimate the second addend in square brackets in right-hand side of (53).
Clearly,

$$
\begin{array}{c}
\int_0^\pi(\nu x^\rho+x^\rho\tau(x))e^{2i\lambda x}G_n(\pi-x)dx=
\int_0^\pi(\nu x^\rho+x^\rho\tau(x))e^{2i\lambda x}Q(\pi-x)\kappa_n(x,\lambda)dx,
\end{array}
$$
where

$$
\begin{array}{c}
\kappa_n(x,\lambda)=e^{2i\lambda(\pi-x)}\int_0^{\pi-x} e^{-2i\lambda t_{3}}P(t_{3})dt_{3}\ldots\\\ldots\int_0^{t_{2n-1}}e^{2i\lambda t_{2n}}Q(t_{2n})dt_{2n}.
\end{array}
$$

It follows from Lemma 1 that
$$
\begin{array}{c}
|\kappa_n(x,\lambda)|\le \int_0^{\pi-x} |P(t_{3})|dt_{3}\ldots\\\ldots\int_0^{t_{2n-1}}|Q(t_{2n})|dt_{2n}\int_0^{t_{2n}} |e^{2i(\pi-x-t_3+\ldots+t_{2n})\lambda }||Q(t_{2n})|dt_{2n}\le\frac{M_2^{2n-2}}{(2n-2)!}.
\end{array}
$$
This and (12) imply
$$
\int_0^\pi(\nu x^\rho+x^\rho\tau(x))e^{2i\lambda x}G_n(\pi-x)dx=\frac{M_3^{2n-1}}{(2n-1)!}\frac{o(1)}{|Im\lambda|^{\rho}}.\eqno(56)
$$
It follows from (53), (56), (55)  that

$$
\sum_{n=1}^\infty|g_{n}(\pi,\lambda)|=\frac{e^{2\pi|Im\lambda|}o(1)}{|Im\lambda|^{\rho}}.\eqno(57)
$$
Relations  (37) and (57) give inequality (45).

{\bf Lemma 5.} {\it Suppose

$$
\lim_{h\to0}\frac{\int_0^hP(x)dx}{h^\rho}=\nu\ne0,\eqno(58)
$$
where $\rho>0$. Then, uniformly in the domain $\Omega_\varepsilon^-$
$$
\begin{array}{c}
\tilde p_{0}(\pi,\lambda)=\frac{\nu\lambda e^{i\pi\lambda}\Gamma(\rho+1)}{2^{\rho}(Re^2\lambda+Im^2\lambda)^{(\rho+1)/2}}e^{i(\rho+1)\arctg(\frac{Re\lambda}{-Im\lambda})}+
\frac{e^{\pi|Im\lambda|}o(1)}{|Im\lambda|^\rho}
\end{array}\eqno(59)
$$
and
$$
\tilde R_{12}(\pi,\lambda)=\frac{e^{\pi|Im\lambda|}o(1)}{|Im\lambda|^\rho}\eqno(60)
$$
as $|Im\mu|\to\infty$.

}
Proof. It follows from (58) that
$$
 \int_0^hP(x)dx=\nu h^\rho+ h^\rho\tau(h),
 $$
 where $\tau(h)\to0$ as $h\to0$.
 Integrating by parts, we obtain
$$
\begin{array}{c}
 \int_0^\pi e^{-2i\lambda t}P(t)dt=e^{-2i\pi\lambda}\int_0^\pi P(y)dy+2i\lambda \int_0^\pi e^{-2i\lambda t}(\int_0^t P(y)dy)dt=\\=
 e^{-2i\pi\lambda}\int_0^\pi P(y)dy+2i\lambda \int_0^\pi e^{-2i\lambda t}[\nu t^\rho+ t^\rho\tau(t)]dt=\\=
 e^{-2i\pi\lambda}\int_0^\pi P(y)dy+2i\nu\lambda \int_0^\infty t^\rho e^{-2i\lambda t}dt-
 2i\nu\lambda\int_\pi^\infty t^\rho e^{-2i\lambda t}dt+2i\lambda \int_0^\pi t^\rho e^{-2i\lambda t}\tau(t)dt.
\end{array}
$$
This together with (10), (11) implies
$$
\int_0^\pi e^{-2i\lambda t}P(t)dt=2i\nu\lambda \int_0^\infty t^\rho e^{-2i\lambda t}dt+\frac{o(1)}{|Im\lambda|^{\rho}}.
$$
Multiplying the last relation by $-ie^{i\lambda \pi}$ and invoking (12), we get (59).

Let us estimate the function $p_0(t,\lambda)$, when $0\le t\le\pi$. Integrating by parts, we obtain

$$
\begin{array}{c}
 \int_0^t e^{-2i\lambda y}P(y)dy=
 e^{-2it\lambda}\int_0^t P(y)dy+2i\lambda \int_0^t e^{-2i\lambda y}dy(\int_0^y P(t_1)dt_1)=\\=
 t^\rho e^{-2it\lambda}(\nu + \tau(t))+2i\lambda \int_0^t y^\rho e^{-2i\lambda y}(\nu + \tau(y))dy.
 \end{array}
$$
This together with (9)  and (13) implies
$$
|p_0(t,\lambda)|\le\frac{c_2}{|Im\lambda|^{\rho}},\eqno(61)
$$
where $c_2$  not depending on $t$.

Let us estimate the function $p_n(t,\lambda)$ if $0\le t\le\pi$.
Since $t_j\ge t_{j+1}$, it follows from (24), (61), and Lemma 1 that
$$
\begin{array}{c}
|p_n(t,\lambda)|=|\int_0^t P(t_1)dt_1\int_0^{t_1}Q(t_2)dt_2\ldots\int_0^{t_{2n-1}}e^{-2i\lambda (t_1-t_2+\ldots+t_{2n-1}-t_{2n})}Q(t_{2n})p_0(t_{2n,\lambda})dt_{2n}|\le\\\le
\frac{c_3}{|Im\lambda|^{\rho}}\int_0^t |P(t_1)|dt_1\int_0^{t_1}|Q(t_2)|dt_2\ldots\int_0^{t_{2n-1}}|Q(t_{2n})|dt_{2n}\le\frac{c_4 M_4^{2n}}{(2n)!|Im\lambda|^{\rho}}.
\end{array}\eqno(62)
$$
It follows from (41), (61), (62), and Lemma 3 that

$$\begin{array}{c}
|p_{n+1}(t,\lambda)|\le\int_0^t|Q(t_2)| |p_n(t_2,\lambda)||e^{2i\lambda t_2}\int_{t_2}^{t}e^{-2i\lambda t_1}P(t_1)dt_1|dt_2
\le\\\le
\max_{0\le t_2\le t}|e^{2i\lambda t_2}\int_{t_2}^{t}e^{-2i\lambda t_1}P(t_1)dt_1|\frac{c_5 M_5^{2n}}{(2n)!|Im\lambda|^{\rho}}=\frac{o(1)M_5^{2n}}{(2n)!|Im\lambda|^{\rho}}
\end{array}\eqno(63)
$$
$(n=0,1\ldots)$, hence,
$$
\begin{array}{c}
\sum_{n=0}^\infty|p_{n+1}(t,\lambda)|=\frac{o(1)}{|Im\lambda|^{\rho}},
\end{array}\eqno(64)
$$
therefore, inequality (60) holds.

Reasoning as above, it is easy to prove the Lemmas 6 and 7. Indeed,
from (22-27), (35), (36) it follows that if in formula (35) for the function $e_{12}(\cdot,\cdot)$ we replace $\lambda$ by  $-\lambda$, swap the functions $P$ and $Q$ and change the sign in front of the function to the opposite, then we get the function $e_{21}(\cdot,\cdot)$. Similarly,
from (28-33), (37), (38) it follows that if in formula (37) for the function $e_{11}(\cdot,\cdot)$ we replace $\lambda$ by  $-\lambda$, swap the functions $P$ and $Q$, then we get the function $e_{22}(\cdot,\cdot)$.

{\bf Lemma 6.} {\it Suppose
$$
\lim_{h\to0}\frac{\int_0^hQ(x)dx}{h^\rho}=\nu\ne0,
$$
where $\rho>0$. Then,
uniformly in the domain $\Omega_\varepsilon^+$
$$
\begin{array}{c}
\tilde q_{0}(\pi,\lambda)=\frac{\nu\lambda e^{-i\pi\lambda}\Gamma(\rho+1)}{2^{\rho}(Re^2\lambda+Im^2\lambda)^{(\rho+1)/2}}e^{i(\rho+1)\arctg(\frac{Re\lambda}{Im\lambda})}+
\frac{e^{\pi|Im\lambda|}o(1)}{|Im\lambda|^\rho}
\end{array}
$$
and
$$
\tilde R_{21}(\pi,\lambda)=\frac{e^{\pi|Im\lambda|}o(1)}{|Im\lambda|^\rho}
$$
as $|Im\mu|\to\infty$.

}

{\bf Lemma 7.} {\it Suppose

$$
\lim_{h\to0}\frac{\int_{\pi-h}^{\pi}Q(x)dx}{h^\rho}=\nu\ne0,
$$
where $\rho>0$. Then,
 uniformly in the domain $\Omega_\varepsilon^-$
$$
\begin{array}{c}
\tilde q_{0}(\pi,\lambda)=\frac{-\nu\lambda e^{i\pi\lambda}\Gamma(\rho+1)}{2^{\rho}(Re^2\lambda+Im^2\lambda)^{(\rho+1)/2}}e^{i(\rho+1)\arctg(\frac{Re\lambda}{-Im\lambda})}+
\frac{e^{\pi|Im\lambda|}o(1)}{|Im\lambda|^\rho},
\end{array}
$$

$$
\tilde R_{12}(\pi,\lambda)=\frac{e^{\pi|Im\lambda|}o(1)}{|Im\lambda|^\rho},
$$
and
$$
e_{22}(\pi,\lambda)=\frac{e^{\pi|Im\lambda|}o(1)}{|Im\lambda|^\rho}
$$
 as $|Im\mu|\to\infty$.
}

Our main result is the following.

{\bf Theorem.} {\it
Suppose
 that at least one of  conditions (65), (66) is satisfied

 $$
 A_{14}\ne0,\eqno(65)
 $$

$$
\lim_{h\to0}\frac{\int_{\pi-h}^{\pi}P(x)dx}{h^{\rho_4}}=\nu_4\ne0,\quad\lim_{h\to0}\frac{\int_0^hQ(x)dx}{h^{\rho_6}}=\nu_6\ne0,
\quad(|A_{13}|+|A_{42}|)|\rho_4-\rho_6|+|-A_{13}\nu_4+A_{42}\nu_6|>0,\eqno(66)
$$
where $\rho_4>0$, $\rho_6>0$;
and also at least one of  conditions (67), (68) is satisfied
$$
A_{32}\ne0,\eqno(67)
$$
$$
\lim_{h\to0}\frac{\int_0^hP(x)dx}{h^{\rho_5}}=\nu_5\ne0,\quad\lim_{h\to0}\frac{\int_{\pi-h}^{\pi}Q(x)dx}{h^{\rho_7}}=\nu_7\ne0,
\quad(|A_{13}|+|A_{42}|)|\rho_5-\rho_7|+|-A_{13}\nu_5+A_{42}\nu_7|>0,\eqno(68)
$$
where $\rho_5>0$, $\rho_7>0$.

Then, the root function system of problem (1), (2) is complete and minimal in $\mathbb{H}$.
}

Proof.
Let $\lambda\in\Omega_\varepsilon^+$ and let $|\lambda|$ be sufficiently large. If (65) holds, then it follows from (5) that

$$
|e_{22}(\pi,\lambda)|\ge c_6e^{\pi|Im \lambda|},\quad|e_{11}(\pi,\lambda)|=o(1)e^{\pi|Im \lambda|},
\quad|e_{12}(\pi,\lambda)|=o(1)e^{\pi|Im \lambda|},\quad|e_{21}(\pi,\lambda)|=o(1)e^{\pi|Im \lambda|}.
$$
This together with (6) implies
$$
|\Delta(\lambda)|\ge c_7e^{\pi|Im\lambda|}.\eqno(69)
$$
 If (66) holds but (65) fails, then
it follows from Lemma 4 that
$$
e_{12}(\pi,\lambda)=\frac{-\nu_4\lambda e^{-i\pi\lambda}\Gamma(\rho_4+1)}{2^{\rho_4}(Re^2\lambda+Im^2\lambda)^{(\rho_4+1)/2}}e^{i(\rho_4+1)\arctg(\frac{Re\lambda}{Im\lambda})}
+\frac{e^{\pi|Im\lambda|}o(1)}{|Im\lambda|^{\rho_4}}.\eqno(70)
$$

$$
e_{11}(\pi,\lambda)=\frac{e^{\pi|Im\lambda|}o(1)}{|Im\lambda|^{\rho_4}}.\eqno(71)
$$
Similarly, by virtue of  Lemma 6, we have
$$
e_{21}(\pi,\lambda)=
\frac{\nu_6\lambda e^{-i\pi\lambda}\Gamma(\rho_6+1)}{2^{\rho_6}(Re^2\lambda+Im^2\lambda)^{(\rho_6+1)/2}}e^{i(\rho_6+1)\arctg(\frac{Re\lambda}{Im\lambda})}+
\frac{e^{\pi|Im\lambda|}o(1)}{|Im\lambda|^{\rho_6}}.\eqno(72)
$$
Relations (6), (66), (70-72) imply
$$
|\Delta(\lambda)|\ge\frac{c_8e^{\pi|Im\lambda|}}{|Im\lambda|^{\max(\rho_4,\rho_6)}}.\eqno(73)
$$

Let $\lambda\in\Omega_\varepsilon^-$ and let $|\lambda|$ be sufficiently large. If (67) holds, then, it follows from (5) that

$$
|e_{11}(\pi,\lambda)|\ge c_9e^{\pi|Im \lambda|},\quad|e_{22}(\pi,\lambda)|=o(1)e^{\pi|Im \lambda|},
\quad|e_{12}(\pi,\lambda)|=o(1)e^{\pi|Im \lambda|},\quad|e_{21}(\pi,\lambda)|=o(1)e^{\pi|Im \lambda|}.
$$
This together with (6) implies
$$
|\Delta(\lambda)|\ge c_{10}e^{\pi|Im\lambda|}.\eqno(74)
$$
 If (68) holds but (67) fails, then it follows from Lemma 5 that
$$
\begin{array}{c}
 e_{12}(\pi,\lambda)=\frac{\nu_5\lambda e^{i\pi\lambda}\nu_5(\rho_5+1)}{2^{\rho_5}(Re^2\lambda+Im^2\lambda)^{(\rho_5+1)/2}}e^{i(\rho_5+1)\arctg(\frac{Re\lambda}{-Im\lambda})}+
\frac{e^{\pi|Im\lambda|}o(1)}{|Im\lambda|^{\rho_5}}.
\end{array}\eqno(75)
$$
Similarly, by virtue of  Lemma 7, we have
$$
\begin{array}{c}
 e_{21}(\pi,\lambda)=\frac{-\nu_7\lambda e^{i\pi\lambda}\Gamma(\rho_7+1)}{2^{\rho_7}(Re^2\lambda+Im^2\lambda)^{(\rho_7+1)/2}}e^{i(\rho_7+1)\arctg(\frac{Re\lambda}{-Im\lambda})}+
\frac{e^{\pi|Im\lambda|}o(1)}{|Im\lambda|^{\rho_7}}.
\end{array}\eqno(76)
$$
and
$$
e_{22}(\pi,\lambda)=\frac{e^{\pi|Im\lambda|}o(1)}{|Im\lambda|^{\rho_7}}.\eqno(77)
$$
Relations (6), (68), (75-77) imply

$$
|\Delta(\lambda)|\ge\frac{c_{11}e^{\pi|Im\lambda|}}{|Im\lambda|^{\max(\rho_5,\rho_7)}}.\eqno(78)
$$
Combining (69), (73), (74), and (78), we obtain
$$
|\Delta(\lambda)|\ge\frac{c_{12}e^{\pi|Im\lambda|}}{|Im\lambda|^{\max(\rho_4,\rho_5,\rho_6,\rho_7)}}
$$
if $\lambda\in\Omega_\varepsilon^-\bigcup\Omega_\varepsilon^+$ ($|\lambda|$ is sufficiently large),
hence, by [2, Th. 2.3] the root function system of problem (1), (2) is complete and minimal in $\mathbb{H}$.

{\bf Remark 2.} If a function $f(x)$ is continuous at the point $0$ or $\pi$, then

$$
\lim_{h\to0}\frac{\int_0^hf(x)dx}{h}=f(0),\quad\lim_{h\to0}\frac{\int_{\pi-h}^\pi f(x)dx}{h}=f(\pi)
$$
correspondingly.
If $f^{(m)}(0)=0, m=0,\ldots n-1$, $f^{(n)}(0)=\nu$ or $f^{(m)}(\pi)=0, m=0,\ldots n-1$, $f^{(n)}(\pi)=\nu$ $(n\ge1)$, then

$$
\lim_{h\to0}\frac{\int_0^hf(x)dx}{h^{n+1}}=\frac{\nu}{(n+1)!},\quad\lim_{h\to0}\frac{\int_{\pi-h}^\pi f(x)dx}{h^{n+1}}=\frac{\nu}{(n+1)!}
$$
correspondingly. Thus, if the functions $P,Q$ are sufficiently smooth and there are certain relations between the values of these functions and their derivatives at the ends of the interval $(0,\pi)$, then the conditions of the Theorem will be satisfied.

\medskip
\medskip
\medskip

\centerline {\bf References}

\medskip
\medskip
\medskip

[1] V.A. Marchenko.  Sturm-Liouville Operators and Their Applications. Birkh\"{a}user, Basel, 1986, Kiev, 1977.

[2] A.A. Lunyov,  M.M. Malamud. On the completeness and Riesz basis property of root subspaces of boundary value problems for first order systems and applications. J. Spectr. Theory. 2015. V. 5. No. 1. P. 17-70.

[3] A.A. Lunyov,  M.M. Malamud. On the completeness of root vectors for first order systems: Application to the Regge problem. Dokl. Math. 2013. V. 88. No. 3. P. 678-683.

[4] M.M. Malamud, L.L. Oridoroga. On the completeness of root subspaces of boundary value problems for first order systems of ordinary differential equations. J. Funct. Anal. 2012. V. 263. No. 7. P. 1939-1980.

[5]  A. V. Agibalova,  M. M. Malamud, L. L. Oridoroga.
On the completeness of general boundary value problems for  $2\times2$
  first order systems of ordinary differential equations.  Methods Funct. Anal. Topol. 2012. V. 18. No. 1. P. 4-18.

[6] A. P. Kosarev,  A. A. Shkalikov.
Spectral asymptotics  of solutions of a $2\times2$
 system of first-order ordinary differential equations.
Math. Notes.  2021. V. 110, No. 6. P. 967-971; translation from Mat. Zametki. 2021. V. 110, No. 6. P. 939-943.

[7] A. Lunyov and M. Malamud. On the Riesz basis property of root vectors system for $2\times2$ Dirac type operators.  J. Math. Anal. Appl. 2016. V. 441. No. 1. P. 57-103.

[8] A.M. Savchuk, I. V. Sadovnichaya. Riesz basis property with parentheses for the Dirac system with summable potential. Sovrem. Math. Fundam. Napravl. 2015. V. 58. P. 128-152.

[9]  A.P. Prudnikov, Yu. A. Brychkov, and O.I. Marichev.  Integraly i ryady (Intagrals and Series), Moscow, Nauka, 1981.

[10] F.G. Tricomi. Integral Equations. 1957. Inc. Publ. NY.

[11] A. M. Savchuk and A. A. Shkalikov. Asymptotic analysis of solutions of ordinary
differential equations with distribution coefficients. Sbornik: Mathematics. 2020. V. 211. No. 11. P. 1623-1659.

\medskip
\medskip
\medskip
\medskip
\medskip

email: alexmakin@yandex.ru

\medskip
\medskip
\medskip
\medskip
\medskip

\end{document}